\documentclass{article}

\usepackage[dvipsnames]{xcolor}

\usepackage{amsmath,amssymb,amsthm}

\usepackage{mathrsfs}

\usepackage{mathtools}
\usepackage{commath}

\usepackage{thmtools}
\usepackage{thm-restate}

\usepackage{cases}
\usepackage{enumitem}
\setlist[enumerate,1]{label=(\roman*)}

\usepackage[pdftex, pdfborderstyle={/S/U/W 0}]{hyperref}
\hypersetup{
    colorlinks=true,
    linkcolor=magenta,
    citecolor=cyan,
}

\usepackage{cleveref}

\usepackage{etoolbox}

\usepackage{comment}

\usepackage[numbers]{natbib}

\numberwithin{equation}{section}

\ifdefined\thmcolor
\declaretheoremstyle[
  shaded={bgcolor=\thmcolor}
]{plain}
\else
\fi

\ifdefined\defcolor
\declaretheoremstyle[
  headfont=\normalfont\bfseries,
  bodyfont=\normalfont,
  shaded={bgcolor=\defcolor}
]{noital}
\else
\declaretheoremstyle[
  headfont=\normalfont\bfseries,
  bodyfont=\normalfont,
]{noital}
\fi

\declaretheorem[style=plain,numberwithin=section,name=Theorem]{theorem}

\declaretheorem[style=plain,sibling=theorem,name=Lemma]{lemma}
\declaretheorem[style=plain,sibling=theorem,name=Corollary]{corollary}
\declaretheorem[style=plain,sibling=theorem,name=Conjecture]{conjecture}
\declaretheorem[style=plain,sibling=theorem,name=Claim]{claim}

\declaretheorem[style=plain,sibling=theorem,name=Question]{question}

\declaretheorem[style=plain,numbered=no,name=Theorem]{theorem-n}
\declaretheorem[style=plain,numbered=no,name=Proposition]{proposition-n}
\declaretheorem[style=plain,numbered=no,name=Lemma]{lemma-n}
\declaretheorem[style=plain,numbered=no,name=Corollary]{corollary-n}
\declaretheorem[style=plain,numbered=no,name=Conjecture]{conjecture-n}
\declaretheorem[style=plain,numbered=no,name=Claim]{claim-n}
\declaretheorem[style=plain,numbered=no,name=Fact]{fact-n}
\declaretheorem[style=plain,numbered=no,name=Open Problem]{openproblem-n}
\declaretheorem[style=plain,numbered=no,name=Question]{question-n}
\declaretheorem[style=plain,numbered=no,name=Observation]{observation-n}

\declaretheorem[style=noital,sibling=theorem,name=Definition]{definition}

\declaretheorem[style=noital,numbered=no,name=Remark]{remark-n}
\declaretheorem[style=noital,numbered=no,name=Definition]{definition-n}
\declaretheorem[style=noital,numbered=no,name=Construction]{construction-n}
\declaretheorem[style=noital,numbered=no,name=Example]{example-n}

\newcommand{\defined}{\mathrel{\coloneqq}}

\undef{\set}
\DeclarePairedDelimiter{\set}{\lbrace}{\rbrace}

\undef{\emptyset}
\newcommand{\emptyset}{\varnothing}

\newcommand{\union}{\mathbin{\cup}}

\newcommand{\from}{\colon}

\undef{\abs}
\DeclarePairedDelimiterX{\abs}[1]
  {\lvert}{\rvert}{\ifblank{#1}{\,\cdot\,}{#1}}

\undef{\norm}
\DeclarePairedDelimiterX{\norm}[1]
  {\lVert}{\rVert}{\ifblank{#1}{\,\cdot\,}{#1}}

\DeclarePairedDelimiterX{\inner}[2]
  {\langle}{\rangle}{\ifblank{#1}{\,\cdot\,}{#1},\ifblank{#2}{\,\cdot\,}{#2}}

\DeclareMathDelimiter{\given}
  {\mathbin}{symbols}{"6A}{largesymbols}{"0C}

\DeclareMathOperator{\Prob}{\mathbb{P}}
\DeclarePairedDelimiterXPP{\prob}[1]
  {\Prob}{\lparen}{\rparen}{}
  {\renewcommand{\given}{\nonscript\;\delimsize\vert\nonscript\;\mathopen{}}#1}

\DeclareMathOperator{\Expec}{\mathbb{E}}
\DeclarePairedDelimiterXPP{\expec}[1]
  {\Expec}{\lparen}{\rparen}{}
  {\renewcommand{\given}{\nonscript\;\delimsize\vert\nonscript\;\mathopen{}}#1}

\DeclareMathOperator{\Var}{Var}
\DeclarePairedDelimiterXPP{\var}[1]
  {\Var}{\lparen}{\rparen}{}
  {\renewcommand{\given}{\nonscript\;\delimsize\vert\nonscript\;\mathopen{}}#1}

\DeclareMathOperator{\Cov}{Cov}
\DeclarePairedDelimiterXPP{\cov}[2]
  {\Cov}{\lparen}{\rparen}{}{#1,#2}

\newcommand{\eps}{\varepsilon}

\newcommand{\EE}{\mathbb{E}}

\usepackage[top=4.34cm]{geometry}

\usepackage{tikz}
\usetikzlibrary {quotes} 
\usepackage{standalone}

\definecolor{comment_green}{RGB}{15, 200, 0}

\begin{document}

\title{Almost colour-balanced spanning forests in complete graphs}

\author{Lawrence Hollom \footnote{\href{mailto:lh569@cam.ac.uk}{lh569@cam.ac.uk}, Department of Pure Mathematics and Mathematical Statistics (DPMMS), University of Cambridge, Wilberforce Road, Cambridge, CB3 0WA, United Kingdom} \and Adva Mond \footnote{\href{mailto:am2759@cam.ac.uk}{am2759@cam.ac.uk}, Department of Pure Mathematics and Mathematical Statistics (DPMMS), University of Cambridge, Wilberforce Road, Cambridge, CB3 0WA, United Kingdom} \and Julien Portier \footnote{\href{mailto:jp899@cam.ac.uk}{jp899@cam.ac.uk}, Department of Pure Mathematics and Mathematical Statistics (DPMMS), University of Cambridge, Wilberforce Road, Cambridge, CB3 0WA, United Kingdom}}



\maketitle

\begin{abstract}
Given $K_n$ whose edges are coloured red and blue, and a forest $F$ of order $n$, we seek embeddings of $F$ with small imbalance, that is, difference between the numbers of red and blue edges.
We show that if the $2$-colouring of the edges of $K_n$ is balanced, meaning that the numbers of red and blue edges are equal, and $F$ has maximum degree $\Delta$, then one can find an embedding of $F$ into $K_n$ whose imbalance is at most $\Delta/2 + 18$, which is essentially best possible and resolves a conjecture of Mohr, Pardey, and Rautenbach.
 Furthermore, we give a tighter bound for the imbalance for small values of $\Delta$.
 In particular, we prove that the imbalance can be taken to be constant in the case where $\Delta<n(1/4 - \eta)$ for any constant $\eta>0$.
\end{abstract}


\section{Introduction}
\label{sec:intro}

Problems concerning substructures appearing when the edges of the complete graph are $2$-coloured have become an area of interest in the past decades.
At one extreme lies Ramsey theory which asks for a monochromatic copy of a given graph $H$. 
The problem of determining $r(H)$ -- the minimal $n$ such that any $2$-colouring of the edges of $K_n$ contains a monochromatic copy of $H$ -- has grown to be a central area of research with remarkable recent breakthroughs by Campos, Griffiths, Morris and Sahasrabudhe~\cite{campos2023exponential} and by Mattheus and Verstraete~\cite{mattheus2024asymptotics}.

Discrepancy is a notion which is closely related to those studied by Ramsey Theory.
Relaxing upon finding a monochromatic copy, one is looking for a copy of a given graph where some colour is much more represented than the other colours.
A significant amount of work has been produced on the subject; see for example \cite{BCJP20, balogh2021discrepancy, erdHos1995discrepancy, GKM22a, GKM22b}.
At the other extreme of Ramsey-type and discrepancy problems, one could ask for a copy of $H$ as far from monochrome as possible, that is to say, a colour-balanced copy of $H$.
We think of the edges of $K_n$ as coloured with colours $-1$ and $+1$, and seek for a copy of $H$ in which the sum of the colours of its edges is close to 0.
It is this question we focus on here, where $H$ is an $n$-vertex forest.

Thus, consider the complete graph $K_n$ with an edge colouring $c : E(K_n) \rightarrow \{-1,+1\}$. 
For convenience we think of the colours blue and red as $-1$ and $+1$, respectively.
For a subgraph $H$ of $K_n$ we define its \emph{sum} to be $c(H) \coloneqq \sum_{e \in E(H)} c(e)$.
We say that $H$ is \emph{balanced} if $c(H) = 0$.
Given a graph $H$, we wish to find a balanced copy of it in $K_n$ equipped with an edge colouring $c : E(K_n) \rightarrow \{-1, +1\}$.
As this needs not be always possible, we seek for such a copy $H'$ of $H$ in $K_n$ such that $|c(H')|$ is as small as possible. For doing so, it makes sense to firstly consider an edge colouring of $K_n$ which is balanced itself, meaning that precisely half of the edges are coloured $-1$ and precisely half of them are coloured $+1$, or in other words, so that $c(K_n)=0$.

In this paper, given a balanced colouring of the edges of $K_n$ by colours $-1$ and $+1$, we consider the problem of embedding an $n$-vertex forest $F$ into $K_n$ such that the colouring induced on $F$ is as balanced as possible.
Some results are known for specific $n$-vertex forests, such as a perfect matching~\cite{ehard2020low, kittipassorn2023existence} (and also~\cite{pardey2022almost} for balanced colouring with $r \ge 3$ colours), factors of $P_3$'s and $P_4$'s~\cite{mohr2022zero} and spanning paths \cite{caro2022zero}.
For the general case where $F$ is any fixed $n$-vertex forest, Mohr, Pardey and Rautenbach posed the following conjecture.

\begin{conjecture}[Conjecture 3 in~\cite{mohr2022zero}]
\label{conj:mohr-spanning-forst}
    Let $c : E(K_n) \rightarrow \{-1,+1\}$ be a balanced colouring and let $F$ be an $n$-vertex forest with maximum degree $\Delta$.
    Then there exists a copy $F'$ of $F$ in $K_n$ with 
    \[\left|c(F') \right| \le \tfrac{1}{2}(\Delta-1).\]
\end{conjecture}

Mohr, Pardey and Rautenbach remarked that one can always find a copy $F'$ of $F$ satisfying $\left|c(F') \right| \le \Delta+1$.
Indeed, swapping the image of any vertex of $F'$ with the image of a leaf induces a change of $c(F')$ by at most $2(\Delta+1)$. 
Therefore, one can `interpolate' between two copies $F_{-}$ and $F_{+}$ of $F$, which satisfy $c(F_{-}) \le 0$ and $c(F_{+}) \ge 0$, by swapping the images of vertices by using a leaf as intermediate.
This way one can obtain a copy $F'$ satisfying $|c(F')| \le \Delta+1$.
See \Cref{lem:trivial-algorithm} for the details of a slight generalisation of this argument.

Pardey and Rautenbach \cite{pardey2023efficiently} made progress towards proving \Cref{conj:mohr-spanning-forst}, showing that for any $\eps >0$ there exists a constant $C_\eps$ such that there is a copy $F'$ of $F$ in any balanced colouring of $E(K_n)$ with colours $-1$ and $+1$ such that $\left|c(F') \right| \le \left(\frac{3}{4} + \eps \right)\Delta +C_\eps$.
In this paper we resolve \Cref{conj:mohr-spanning-forst} up to an additive constant term.

\begin{theorem}
\label{thm:spanning-forest-general}
    Let $F$ be an $n$-vertex forest with maximum degree $\Delta$, and let $c : E(K_n) \rightarrow \{-1, +1\}$ be a balanced colouring of the edges of $K_n$. Then there exists a copy $F'$ of $F$ in $K_n$ such that
    \begin{align}
    \label{eq:conjUB}
        \left|c(F') \right| \le \tfrac{1}{2}\Delta + 18.
    \end{align}
\end{theorem}

In fact, we prove the following stronger result for certain values of $\Delta$, from which we derive \Cref{thm:spanning-forest-general}.

\begin{theorem}
\label{thm:spanning-forest}
    Let $F$ be an $n$-vertex forest with maximum degree $\Delta$, and let $c : E(K_n) \rightarrow \{-1, +1\}$ be a balanced colouring of the edges of $K_n$.
    Then there exists a copy $F'$ of $F$ in $K_n$ such that
    \begin{align}
    \label{eq:spanning-forest}
        \left|c(F') \right|\leq \begin{cases}
            \frac{1}{2}\Delta + 9 & \text{ if } \Delta \ge \frac{1}{2}n \text{ or } \Delta \le 15, \\
            \frac{1}{2}\left(\Delta - \frac{1}{4}n + 3 \right) + \sqrt{\frac{1}{4}\left(\Delta - \frac{1}{4}n + 1 \right)^2 + 4n} & \text{ if } 15 < \Delta < \frac{1}{2}n.
        \end{cases}
    \end{align}
\end{theorem}

For values of $\Delta$ almost all the way up to $\frac{n}{4}$, \Cref{thm:spanning-forest} gives a constant upper bound.
\begin{corollary}
\label{cor:main}
    Let $F$ be an $n$-vertex forest with maximum degree $15 \le \Delta \le \left(\frac{1}{4} - \eta \right)n$, for some constant $\eta \in \left(0, \frac{1}{4} \right)$.
    Let $c : E(K_n) \rightarrow \{-1, +1\}$ be a balanced colouring of the edges of $K_n$.
    Then there exists a copy $F'$ of $F$ in $K_n$ such that
    \begin{align*}
        \left|c(F') \right| \le (1+o(1))\frac{4}{\eta}.
    \end{align*}
\end{corollary}

While the additive constant of $-\frac{1}{2}$ in \Cref{conj:mohr-spanning-forst} cannot hold in general, as the example of embedding a Hamiltonian path into $K_n$ where $n$ is even shows (the absolute imbalance would have to be at least $1$ due to parity considerations), we have made little effort to optimise the additive constant $18$ in our bound.

We conclude this introduction by mentioning related work in the field.
Previous research has established the existence of spanning forests with small imbalance, although without specifying the isomorphic class of the forest.
For instance, Caro and Yuster~\cite{caro2016zero} studied edge-colourings of $K_n$ with colours $-1$ and $+1$ which have small imbalance, and show tight bounds on the imbalance that allows for the existence of an almost zero-sum spanning tree $T$, but without fixing the isomorphic class of $T$.
Caro, Hansberg, Lauri and Zarb~\cite{caro2022zero} extended this study to sparser graphs $G$ in place of the complete graph $K_n$.
Finally, F\"{u}redi and Kleitman~\cite{furedi1992zero} showed that when the edges of $K_n$ are coloured using the set $\{0, \dots, n-1\}$, the graph contains a spanning tree whose sum is $0$ mod $n-1$, a result for which Schrijver and Seymour later provided a simpler proof~\cite{schrijver1991simpler}.




\subsection{Notation and terminology}
\label{subsec:notation}

Throughout our proofs, we will often need to refer to the number of edges of a particular colour or within a particular subgraph, and to this end we use the following notation.
We associate the colour blue with $-1$ and the colour red with $+1$, and we will, without comment, switch between conventions as convenient.

For a graph $H$ and a vertex $v \in V(H)$ we denote by $N_H(v)$ the neighbourhood of $v$ in $H$ and by $\deg_H(v) = \left|N_H(v) \right|$ the degree of $v$ in $H$.
Similarly, for a subgraph $H' \subset H$ we respectively denote by $N_{H'}(v)$ and $\deg_{H'}(v) = \left|N_{H'}(v) \right|$ its neighbourhood and degree in $H'$.
Given a colouring $c : E(K_n) \rightarrow \{-1,+1\}$, we denote by $B$ and $R$ the subgraphs of $K_n$ spanned by the edges coloured $-1$ (blue) and $+1$ (red), respectively.
In particular, $N_B(v)$ and $N_R(v)$ stand for the blue and red neighbourhoods and $\deg_B(v)$ and $\deg_R(v)$ stand for the blue and red degrees of $v$ in $K_n$ under $c$, respectively.

As mentioned earlier, we say that a colouring $c : E(K_n) \rightarrow \{-1,+1\}$ is balanced if we have $c(K_n) = 0$, where we recall that for a graph $H$ its sum $c(H)$ is defined to be $c(H):=\sum_{e\in E(H)}c(e)$.
We will slightly abuse notation and say that a graph $H$ is \emph{(almost) colour-balanced} or simply \emph{(almost) balanced} if the induced colouring on $H$ contains (almost) as many blue edges than red edges, where the term ``almost'' is not defined quantitatively and serves mostly as a heuristic.
Given a graph $H$ we denote by $\Delta(H)$ and $\delta(H)$ its maximum and minimum degrees.
For a real number $x \in \mathbb R$ we write $(x)_+ = \max(x,0)$.


\section{Preliminaries}
\label{sec:preliminaries}

We start with stating and proving \Cref{lem:trivial-algorithm}.
This lemma allows us to find an almost balanced embedding, provided that we already know of the existence of two embeddings: one with a non-negative sum, and one with a non-positive sum.
This is a slight strengthening of a lemma appearing in previous papers (see, e.g., \cite{mohr2022zero}).

\begin{lemma}
\label{lem:trivial-algorithm}
    Let $n$ be an integer and suppose that $c ~:~ E(K_n) \rightarrow \{-1,+1\}$ is a balanced colouring.
    Let $H$ be a graph on $n$ vertices with minimum degree $\delta_H \coloneqq \delta(H)$.
    Let $h_B, h_R ~:~ V(H) \rightarrow V(K_n)$ be two embeddings of $H$ into $K_n$, with images $H_B, H_R$, such that $c(H_B) \le 0 \le c(H_R)$.
    Let $I \coloneqq \left\{v \in V(H) ~:~ h_B(v) \neq h_R(v) \right\}$ and set $\Delta_I = \max\left\{\deg_H(v) ~:~ v \in I \right\}$.
    Then there exists a subgraph $H'$ of $K_n$ which is isomorphic to $H$ and such that $|c(H')| \le \Delta_I + \delta_H$.
\end{lemma}
\begin{proof}
    Since $c(H_B) \le 0 \le c(H_R)$, we may assume that $c(H_R) > \Delta_I+\delta_H$ and $c(H_B) < -(\Delta_I+\delta_H)$, as otherwise we are done by choosing $H'$ to be either $H_R$ or $H_B$.
    We define a sequence of embeddings $h_R = h_0, h_1, \ldots, h_{\ell} = h_B$ of $H$ into $K_n$, for some integer $\ell \ge 1$, such that for each $i \in [\ell]$, $h_i$ is obtained from $h_{i-1}$ by switching the images of precisely two vertices in $I$.
    Note that this is possible, as there is a permutation $\sigma : V(H) \rightarrow V(H)$ of the vertices of $H$ for which $h_R(\sigma(v)) = h_B(v)$ for all $v \in V(H)$.
    More precisely, we define $\sigma$ to be a product of swaps, all of which fix all elements in $V(H)$, including $I$, except for two vertices in $I$ which are swapped.
    Let $\sigma_1, \ldots, \sigma_{\ell}$ be this sequence of swaps, and for each $i \in [\ell]$ let $u_i, v_i \in I$ be such that $\sigma_i(u_i) = v_i$, $\sigma_i(v_i) = u_i$ and everything else is fixed.
    
    Let $i \in [\ell]$.
    Given $h_{i-1}$, we define $h_i$ as follows.
    Assume first that, without loss of generality, we have $\deg_H(v_i) = \delta_H$ (the case where $\deg_H(u_i) = \delta_H$ is similar).
    Then we set
    \begin{align*}
        h_i(x) \defined \begin{cases}
            h_{i-1}(v_i) & \text{ if } x=u_i, \\
            h_{i-1}(u_i) & \text{ if } x=v_i, \\
            h_{i-1}(x) & \text{ otherwise,}
        \end{cases}
    \end{align*}
    i.e., we swap the images of $u_i$ and $v_i$.
    Let $H_i$ be the image of the graph $H$ under $h_i$, for all $i \in [\ell]$.
    Note that in this case we have
    \begin{align*}
        \left|c(H_i) - c(H_{i-1}) \right| \le 2(\deg_H(u_i) + \deg_H(v_i)) \le 2(\Delta_I + \delta_H).
    \end{align*}
    Assume now that $\deg_H(v_i) > \delta_H$ and $\deg_H(u_i) > \delta_H$.
    In this case we swap the images of $u_i$ and $v_i$ in three steps using two auxiliary functions.
    Let $w \in V(H) \setminus \{u_i, v_i\}$ be a vertex with $\deg_H(w) = \delta_H$.
    We define
    \begin{align*}
        h'_i(x) \defined \begin{cases}
            h_{i-1}(u_i) & \text{ if } x=w, \\
            h_{i-1}(w) & \text{ if } x=u_i, \\
            h_{i-1}(x) & \text{ otherwise,}
        \end{cases}
    \end{align*}
    \begin{align*}
        h''_i(x) \defined \begin{cases}
            h'_i(v_i) & \text{ if } x=w, \\
            h'_i(w) & \text{ if } x=v_i, \\
            h'_i(x) & \text{ otherwise,}
        \end{cases}
    \end{align*}
    and
    \begin{align*}
        h_i(x) \defined \begin{cases}
            h''_i(u_i) & \text{ if } x=w, \\
            h''_i(w) & \text{ if } x=u_i, \\
            h''_i(x) & \text{ otherwise.}
        \end{cases}
    \end{align*}
    That is, $h'_i$ is $h_{i-1}$ with the images of $u_i$ and $w$ swapped, $h''_i$ is $h'_i$ with the images of $v_i$ and $w$ swapped and, finally, $h_i$ is $h''_i$ with the images of $u_i$ and $w$ swapped.
    Note that $h_i$ agrees with $h_{i-1}$ everywhere except for swapping $u_i$ and $v_i$, and in particular $h_i(w) = h_{i-1}(w)$.
    Moreover, if we denote by $H_i, H'_i$ and $H''_i$ the images of $H$ under $h_i, h'_i$ and $h''_i$, respectively, we have
    \begin{align*}
        \left|c(H'_i) - c(H_{i-1}) \right| \le 2(\deg_H(u_i) + \deg_H(w)) \le 2(\Delta_I + \delta_H),
    \end{align*}
    \begin{align*}
        \left|c(H''_i) - c(H'_i) \right| \le 2(\deg_H(w) + \deg_H(v_i)) \le 2(\Delta_I + \delta_H),
    \end{align*}
    and
    \begin{align*}
        \left|c(H_i) - c(H''_i) \right| \le 2(\deg_H(u_i) + \deg_H(w)) \le 2(\Delta_I + \delta_H).
    \end{align*}
    
    In any case for every $i \in [\ell]$ we have $h_i(\sigma_i(v)) = h_{i-1}(v)$ for all $v \in V(H)$, and therefore $h_0 = h_R$ and $h_{\ell} = h_B$.
    As $c(H_R) > \Delta_I + \delta_H$ and $c(H_B) < -(\Delta_I + \delta_H)$, by the above we get that there exist $i \in [\ell]$ for which $|c(\hat{H})| \le \Delta_I + \delta_H$, where $\hat{H} \in \left\{H_i, H'_i, H''_i \right\}$.
\end{proof}

Consider an embedding $f : V(F) \rightarrow V(K_n)$.
We introduce the following definition.

\begin{definition}
\label{def:sgnfix}
The \emph{sign} of $f$ is defined to be the sign of $c(f) = \sum_{e \in E(F)} c(f(e))$.
Furthermore, given a subset $U \subset V(K_n)$ we say that a subset $L \subset V(F)$ is \emph{$U$-sign-fixing for $F$} if any two embeddings $f, f' : V(F) \rightarrow V(K_n)$ such that $f(L), f'(L) \subset U$ have the same sign whenever $f|_L = f'|_L$.
\end{definition}
Note that, given a subset $U \subset V(K_n)$, if a subset $L \subset V(F)$ is $U$-sign-fixing for $F$, then so is any subset $L' \subset V(F)$ with $|L'| \le |U|$ which contains $L$.
In our proof of \Cref{thm:spanning-forest} we consider a certain subset $U \subset V(K_n)$ and we will show that if a set $L$ of vertices of large degrees is \emph{not} $U$-sign-fixing, then we can find an embedding of $F$ with a small sum (in absolute value) using \Cref{lem:trivial-algorithm}.
We then show that if $L$ \emph{is} $U$-sign-fixing, then, by considering a minimal $U$-sign-fixing set $M \subset L$, we can derive the desired bound.

Once a $U$-sign-fixing set $L$ of vertices is embedded in $U$, the sign of any embedding of $F$ which extends $L$ is already determined, no matter how the remaining vertices are embedded.
In such a situation, we will also refer to the sign of an embedding $g : L \rightarrow U$, and we will say that $g$ is \emph{red} if any embedding $g'$ of $F$ extending $g$ satisfies $c(g'(F)) \ge 0$, and we say that $g$ is \emph{blue} if $c(g'(F)) \le 0$ for any such embedding $g'$.

\begin{definition}
\label{def:balancedvtx}
Given a $2$-colouring of the edges of $K_n$ and an integer $r$, we say that a vertex in $K_n$ is \emph{$r$-balanced} if it is adjacent to at least $r$ edges in each colour.
\end{definition}
In the following lemma we show that if $r$ is not too big, one can find many $r$-balanced vertices in a balanced colouring of $K_n$.
\begin{lemma}
\label{lem:ExistenceCbalancedVertices}
Let $n \ge 5$ be an integer and $\frac{1}{n} \le \eps < \frac{1}{4}$.
Then in any balanced $2$-colouring of the edges of $K_n$ by red and blue there are at least $\varepsilon n + 1$ many vertices which are $\left(\frac{1}{4} - \varepsilon \right)n$-balanced.
\end{lemma}
We remark that a version of \Cref{lem:ExistenceCbalancedVertices} for $\eps=0$ is Theorem 2 of \cite{mohr2022zero}.
The proof of \Cref{lem:ExistenceCbalancedVertices} is somewhat technical and not relevant to our other arguments, and so we include it in \Cref{appendixA}.


\section{Embedding an $n$-vertex forest}
\label{sec:spanning-forest}

We prove \Cref{thm:spanning-forest} by dealing separately with the cases $\Delta \le 15$, $15 < \Delta \le \frac{1}{2}n$ and $\Delta > \frac{1}{2}n$.
In the proof of the second case, the following lemma will be the crux of our argument.
Recall that for a real number $x \in \mathbb R$ we denote $(x)_+ = \max(x,0)$.

\begin{lemma}
\label{lem:Delta<=n/2}
    Let $c : E(K_n)\rightarrow \{-1,+1\}$ be a balanced colouring of the edges of $K_n$, for some integer $n \ge 32$, and let $F$ be an $n$-vertex forest with maximum degree $\Delta \le \frac{1}{2}n$.
    Let $\frac{1}{n} \le \eps \le \frac{1}{8}$.
    Then there is a copy $F'$ of $F$ in $K_n$ such that
    \begin{align}
    \label{eq:Delta<=n/2}
        |c(F')| \leq \max\left(1+\tfrac{2}{\eps}, 2 + \left( \Delta - \left(\tfrac{1}{4} - 2\eps \right)n\right)_+\right),
    \end{align}
\end{lemma}

\begin{proof}
    By \Cref{lem:ExistenceCbalancedVertices} we know that there exists a subset $U \subset V(K_n)$ of at least $\eps n + 1$ many $\left(\frac{1}{4}-\eps \right)n$-balanced vertices, under the colouring $c$.
    Define the set $L_{\eps} \subset V(F)$ of vertices of `large degree' to be
    \[L_{\eps} \coloneqq \left\{v \in V(F) ~:~ \deg_F(v) \ge \tfrac{2}{\eps} \right\}. \]
    Note that $|L_{\eps}| \le \frac{2n}{2/\eps} = \eps n$.
    We separate the proof into two cases.
    If $L_{\eps}$ is not $U$-sign-fixing then we find a copy $F'$ of $F$ in $K_n$ with $|c(F')| \le 1 + \frac{2}{\eps}$ by using \Cref{lem:trivial-algorithm}.
    On the other hand, if $L_{\eps}$ is $U$-sign-fixing then we consider a minimal $U$-sign-fixing set $M \subset L_{\eps}$ and we find a copy $F'$ of $F$ with $|c(F')| \le 2 + (\Delta - (\tfrac{1}{4} - 2\eps )n)_+$ using the definition of $U$.\\

    Suppose first that $L_{\eps}$ is not $U$-sign-fixing.
    Then there are embeddings $h_B, h_R$ of $F$ with $h_B|_{L_{\eps}} = h_R|_{L_{\eps}}$ and $h_B(L_{\eps}), h_R(L_{\eps}) \subset U$ such that $c(h_B(F)) \le 0 \le c(h_R(F))$.
    Note further that, by definition, we have
    \[I \coloneqq \left\{v \in V(F) ~:~ h_R(v) \neq h_B(v) \right\} \subset V(F)\setminus L_{\eps}, \]
    implying that $\Delta_I \coloneqq \max\{\deg_F(v) ~:~ v \in I \} \le \frac{2}{\eps}$. 
    Since $\delta_F = 1$, by \Cref{lem:trivial-algorithm}, we find that there is an embedding $h$ of $F$ with $F' \coloneqq h(F)$ for which
    \begin{align}
    \label{eq:not-sign-fixing-bound}
        \abs{c(F')}\le \tfrac{2}{\eps}+1,
    \end{align}
    so in particular (\ref{eq:Delta<=n/2}) holds in this case.

    Suppose now that $L_{\eps}$ is $U$-sign-fixing for $F$, and let $M \subset L_{\eps}$ be a minimal $U$-sign-fixing set.
    Consider the subset
    \begin{align}
    \label{eq:definition-of-n}
        N \coloneqq \left\{v \in M ~:~ \deg_{F[M]}(v) \ge 2 \right\}.
    \end{align}
    As $c$ is a balanced colouring we have that $M$ is non-empty, and moreover, as $F$ is a forest we have $N \subsetneq M$.
    Moreover, since $M$ is a minimal $U$-sign-fixing set for $F$, and since $F$ is a forest, we know that $N$ is not $U$-sign-fixing.
    Hence, there exist $f, g : M \rightarrow U$ with $f|_N = g|_N$ and such that $f$ is red and $g$ is blue.
    That is, for every $f', g' : V(F) \rightarrow V(K_n)$ such that $f'|_M = f$ and $g'|_M = g$ we have $c(f') \ge 0$ and $c(g') \le 0$.
    Note that we have $|M| \le |L_{\eps}| \le \eps n$ and recall that $|U| \ge \eps n + 1$.
    In particular, there exists a $\left(\frac{1}{4} - \eps \right)n$-balanced vertex $a \in U \setminus g(M)$.

    Label the vertices of $M \setminus N$ as $\{v_1, \ldots, v_r \}$, where $r \ge 1$ is some integer.
    As $a \notin g(M)$ we also have $a \notin g(N) = f(N)$.
    However, it is possible to have $a \in f(M\setminus N)$.
    If this is the case, then we may assume, without loss of generality, that $f(v_r) = a$.
    
    We now define a sequence of embeddings $h_0, h_1, \ldots, h_{3r}$ of $M$ into $U$.
    Intuitively, we define these embeddings to `interpolate' between $g$ and $f$, with $h_{k-1}$ and $h_k$ differing only in where they send a single vertex $v$.
    The technicalities arise due to the possibility that the images of $f$ and $g$ may overlap, and so we must ensure that each $h_k$ is indeed an injection.

    Define the functions $h_k :  M \rightarrow V(K_n)$ inductively, as follows.
    Firstly set $h_0(v_j) \coloneqq g(v_j)$ for all $j \in [r]$, and for all $k \in [3r]$ set $h_k|_N = f|_N = g|_N$.
    Then for each $i \in [r]$ in turn define
    \begin{align*}
    h_{3i-2}(v_j) &\defined \begin{cases}
        a & \text{ if } j > i, h_{3i-3}(v_j)=f(v_i)\neq h_{3i-3}(v_i), \\
        h_{3i-3}(v_j) & \text{ otherwise},
    \end{cases}\\
    h_{3i-1}(v_j) &\defined \begin{cases}
        f(v_i) & \text{ if } j = i, \\
        h_{3i-2}(v_j) & \text{ otherwise},
    \end{cases}\\
    h_{3i}(v_j) &\defined \begin{cases}
        h_{3i-3}(v_i) & \text{ if } j > i, h_{3i-1}(v_j)=a, \\
        h_{3i-1}(v_j) & \text{ otherwise}.
    \end{cases}
    \end{align*}
    
    In words, in the $i$th step, these three functions ensure that $v_i$ is mapped to $f(v_i)$. 
    As we wish to `interpolate' between $g$ and $f$ using injections only, we must make sure that no other vertex except for $v_i$ is sent to $f(v_i)$.
    If this is not the case, and some vertex $v_j$ is already mapped to $f(v_i)$ at this point, then $h_{3i-2}$ moves it to $a$, where it is temporarily sent.
    Then $h_{3i-1}$ sends $v_i$ to $f(v_i)$, and $h_{3i}$ sends $v_j$ to where $v_i$ used to be sent, freeing up the vertex $a$.
    Note that if $f(v_i)\notin h_{3i-3}(M)$, then $h_{3i-2}=h_{3i-3}$ and $h_{3i}=h_{3i-1}$.
    We now prove that these embeddings $h_i$ have the properties we desire.

    \begin{claim}
    \label{claim:interpolating-embeddings-are-good}
        For all $k \in [3r]$ we have $h_k(M) \subset U$, where $h_k$ is injective, and the functions $h_{k-1}$ and $h_k$ differ in where they embed at most one vertex.
        Furthermore, for all $1\leq i < r$, $a\notin h_{3i}(M)$, and for all $1\le j \le i\le r$, $h_{3i}(v_j)=f(v_j)$.
    \end{claim}
    
    \begin{proof}
        It is immediate from the definition that $h_k$ is injective, that $h_{k-1}$ and $h_k$ differ in where they embed at most one vertex.
        Moreover, as $g(M), f(M) \subset U$ and since for all $k \in [3r]$ and all $v \in M$ we have $h_k(v) \in \{g(v), f(v), a \}$, we get that $h_k(M) \subset U$ as well.
        
        We prove the rest of the claim inductively, noting that $h_0$ does indeed satisfy all of the required properties.
        Let $i\in [r]$ and assume that $a \notin h_{3i-3}(M)$ and $h_{3i-3}(v_j) = f(v_j)$ for all $j \le i-1$.
        We prove that the same holds for $h_{3i}$ as well.
        As, by their definition, $h_{3i}(v_j) = h_{3i-3}(v_j) = f(v_j)$ for all $j \le i-1$, it is enough to show that $h_{3i}(v_i) = f(v_i)$ and that if $i < r$ then $a \notin h_{3i}(M)$.
        
        If we already have $h_{3i-3}(v_i) = f(v_i)$ then we get $h_{3i-3} = h_{3i-2} = h_{3i-1} = h_{3i}$ and we are done, so we may assume that this is not the case.
        If $f(v_i) \notin h_{3i-3}(M)$ then we get $h_{3i-3} = h_{3i-2}$, and then $h_{3i-1}(v_i) = f(v_i)$.
        As we further have in this case that $h_{3i} = h_{3i-1}$ and $a \notin h_{3i-1}(M)$, we are done.

        Hence, we may assume that $f(v_i) \in h_{3i-3}(M)$ and that $h_{3i-3}(v_j) = f(v_i)$ for some $j > i$.
        In this case we have $h_{3i-1}(v_j) = h_{3i-2}(v_j)$ for all $j\neq i$.
        By definition we have $h_{3i-1}(v_i) = f(v_i)$, and moreover, we have $h_{3i}(v_j) = h_{3i-1}(v_j)$ for all $j \le i$, so in particular $h_{3i}(v_i) = f(v_i)$.
        Further, note that if $i=r$ then we may have $f(v_r) = a$ in which case $a \in h_{3r}(M)$.
        If $i < r$, then although $a \not\in h_{3i-3}(M)$, we may still have $a \in h_{3i-2}(M)$ or $a \in h_{3i-1}(M)$.
        This happens precisely when there exists $j > i$ for which $h_{3i-3}(v_j) = f(v_i)$.
        However, we then get that $h_{3i-1}(v_j) = h_{3i-2}(v_j) = a$, but then $h_{3i}(v_j) = h_{3i-3}(v_j) \neq a$, as $a \notin h_{3i-3}(M)$, so in particular $a \notin h_{3i}(M)$.
        This completes the induction, and hence the proof of \Cref{claim:interpolating-embeddings-are-good}.
    \end{proof}

    Recall that $M$ is sign fixing, so each $h_i$ is either blue or red, and we know that $h_0 = g$ is blue and $h_{3r} = f$ is red.
    Hence, there exists an index $t \in [3r]$ for which $h_{t-1}$ is blue and $h_t$ is red.
    By \Cref{claim:interpolating-embeddings-are-good}, $h_{t-1}$ and $h_t$ differ only in the image of one vertex $v \in M\setminus N$, as $h_{t-1}|_N = h_t|_N$, so in particular, recalling \eqref{eq:definition-of-n}, we have $\deg_M(v) \le 1$.
    Moreover, by \Cref{claim:interpolating-embeddings-are-good} we know that $h_{t-1}(v) \in U$, meaning it is $\left(\frac{1}{4} - \eps \right)n$-balanced.

    We now define embeddings $h'_{t-1}$ and $h'_t$ of $F$ into $K_n$, extending $h_{t-1}$ and $h_t$.
    Our goal is to do so in such a way that the difference $c(h'_t(F)) - c(h'_{t-1}(F))$ is bounded from above.
    As $h_{t-1}$ is blue and $h_t$ is red, we are guaranteed to have $c(h'_{t-1}(F)) \le 0 \le c(h'_t(F))$ no matter how we choose these two extensions.
    We consider the vertex $h'_{t-1}(v)$, which we know to be a $\left(\frac{1}{4} - \eps \right)n$-balanced vertex.
    We then define an extension $h_{t-1}$ so that, locally, it contains as many red edges as possible in the red neighbourhood of $h'_{t-1}(v)$.
    As $h'_{t-1}(v)$ is $\left(\frac{1}{4}-\eps \right)n$-balanced, the number of red neighbours of $h'_{t-1}(v)$ in $K_n$ outside of $h_{t-1}(M) \cup \{h_t(v)\} \setminus \{h_{t-1}(v) \}$ is at least
    \begin{align}
    \label{eq:Rnbrs}
        \left(\tfrac{1}{4}-\eps \right)n - |M| \ge \left(\tfrac{1}{4} - 2\eps \right)n.
    \end{align}
    Let $h'_{t-1} : V(F) \rightarrow V(K_n)$ be an embedding of $F$ such that $h'_{t-1}|_M = h_{t-1}$ and with the following two properties.
    The first is that $h'_{t-1}$ has as many red neighbours of $h_{t-1}(v)$ outside $h_{t-1}(M)$ as possible; more formally, with either $h'_{t-1}(N_F(v) \setminus M) \subset N_R(h_{t-1}(v))$ or $N_R(h_{t-1}(v)) \subset h'_{t-1}(N_F(v) \setminus M)$.
    Thus, by (\ref{eq:Rnbrs}), we get that the number of blue edges in $h'_{t-1}(F) \subset K_n$ from $h_{t-1}(v)$ to a vertex not in $h_{t-1}(M) \cup \{h_t(v)\} \setminus \{h_{t-1}(v) \}$ is equal to
    \begin{align}
    \label{eq:num-blue-nbrs-outside-of-M}
        e_B\left(h_{t-1}(v), h'_{t-1}(N_F(v)) \setminus \left(h_{t-1}(M) \cup \{h_t(v) \} \right) \right) \le \left(\Delta - \left(\tfrac{1}{4}-2\eps \right)n\right)_+.
    \end{align}
    The second property that $h'_{t-1}$ has to satisfy is that it maps a leaf of $F$ to $h_t(v)$.
    Note that, as $h_{t-1}$ and $h_t$ differ only at where they map $v$ to, we have $h_t(v) \notin h_{t-1}(M)$.
    Hence, we may indeed map a leaf of $F$ to $h_t(v)$ under $h'_{t-1}$. 
    We denote $w \coloneqq (h'_{t-1})^{-1}(h_t(v))$.

    Having defined $h'_{t-1}$, we define $h'_t$ to be an extension of $h_t$, as follows.
    \begin{align*}
        h'_t(x) \coloneqq \begin{cases}
            f(v) & \text{ if } x=v, \\
            h'_{t-1}(v) & \text{ if } x=w, \\
            h'_{t-1}(x) & \text{ otherwise,}
        \end{cases}
    \end{align*}
    i.e. $h'_t$ is equal to $h'_{t-1}$ but with $v$ and $w$ swapped.
    Since $h_{t-1}$ is blue and $h'_{t-1}|_M = h_{t-1}$ we know that $c(h'_{t-1}(F)) \le 0$.
    Similarly, as $h_t$ is red and $h'_t|_M = h_t$ we know that $c(h'_t(F)) \ge 0$.
    However, by (\ref{eq:num-blue-nbrs-outside-of-M}) we know that all but at most $\left(\Delta - \left(\frac{1}{4} - 2\eps \right)n \right)_+$ many edges from $v$ to outside of $M$ were embedded by $h'_{t-1}$ to red edges.
    Thus, we get that
    \begin{align}
    \label{eq:diff-h't-h't-1}
        c(h'_t(F)) - c(h'_{t-1}(F)) \le 2\deg_M(v) + 2\deg_F(w) + 2\left(\Delta - \left(\tfrac{1}{4} - 2\eps \right)n \right)_+.
    \end{align}
    Recall that $w$ is a leaf of $F$ and that $v \in M\setminus N$ so we have $\deg_M(v) \le 1$.
    Hence, by (\ref{eq:diff-h't-h't-1}), there exists an embedding $h$ of $F$ with $F' \coloneqq h(F)$ (which is either $h'_{t-1}$ or $h'_t$), for which we have
    \begin{align}
    \label{eq:diff}
        |c(F')| \le 2 + \left(\Delta - \left(\tfrac{1}{4} - 2\eps \right)n \right)_+,
    \end{align}
    finishing the proof.
\end{proof}

For the proof of \Cref{thm:spanning-forest} we need one more technical lemma to deal with the case $15 < \Delta < \frac{1}{2}n$.
More precisely, \Cref{lem:Delta<=n/2} implies \Cref{thm:spanning-forest} when $15 < \Delta < \frac{1}{2}n$ by optimising (\ref{lem:Delta<=n/2}) over $\eps$, as given in the following lemma.
\begin{lemma}
\label{lem:optim}
    Let $\xi \in \left[-\frac{1}{4} + \frac{16}{n}, \frac{1}{4} \right]$ and $n \ge 32$ be an integer.
    Then
    \begin{align}
    \label{eq:optim}
        \phi_n(\xi) &\coloneqq \inf_{\frac{1}{n} \le \eps \le \frac{1}{8}} \max \left(\tfrac{2}{\eps}+1, 2 + (2\eps n + \xi n)_+ \right) \le
        \frac{3+\xi n}{2} + \frac{\sqrt{(1+\xi n)^2+16n}}{2}.
    \end{align}
\end{lemma}
The proof of \Cref{lem:optim} is purely technical and consists only of elementary techniques.
As it does not serve to enlighten the rest of the paper, it is given in \Cref{appendixB}.

We are now ready to prove our main result.
As mentioned earlier, the case $\Delta \le \frac{1}{2}n$ follows immediately from \Cref{lem:Delta<=n/2,lem:optim}.
The case $\Delta \ge \frac{1}{2}n$ is obtained by explicitly embedding up to three vertices of $F$ of largest degrees.

\begin{proof}[Proof of \Cref{thm:spanning-forest}]
Firstly, note that if $\Delta \le 15$ then \Cref{lem:trivial-algorithm} immediately implies (\ref{eq:spanning-forest}). 
Hence, we assume that $\Delta > 15$ and we continue the proof by further splitting into two cases according to the value of $\Delta$.

\paragraph*{Case 1: $15 < \Delta < \frac{1}{2}n$.}
Note that $n \ge 32$.
Let $\xi \in \left[-\frac{1}{4} + \frac{16}{n}, \frac{1}{4} \right]$ be such that $\Delta = \left(\frac{1}{4} + \xi \right)n$.
Given $\frac{1}{n} \le \eps \le \frac{1}{8}$, we know by \Cref{lem:Delta<=n/2} that there exists a copy $F'$ of $F$ in $K_n$ for which
\[\left|c(F') \right| \le \max \left(1+\tfrac{2}{\eps}, 2 + (2\eps n + \xi n)_+ \right). \]
As $\xi n = \Delta - \frac{1}{4}n$, by optimising over all possible values of $\eps$, \Cref{lem:optim} implies (\ref{eq:spanning-forest}) for this case.

\paragraph*{Case 2: $\Delta \ge \frac{1}{2}n$.}
We explicitly embed a few vertices of largest degrees in $F$.
Let $v_1, v_2, \ldots, v_n$ be an ordering of the vertices of $F$ such that $\Delta = \deg_F(v_1) \ge \deg_F(v_2) \ge \cdots \ge \deg_F(v_n)$, and set $\Delta' = \deg_F(v_2)$.
Let $\eps = \frac{1}{n}$.
By \Cref{lem:ExistenceCbalancedVertices} we know that there exists a vertex $x \in V(K_n)$ which is $\left(\frac{1}{4} - \eps \right)n$-balanced.
We proceed with simple analyses of two subcases.

\paragraph*{Case 2.1: $\Delta' < \frac{1}{2}\Delta$.}
\label{case:Delta-prime-small}
We embed $v_1$ at $x$ and consider an extension of this embedding $f : V(F) \rightarrow V(K_n)$ chosen uniformly at random.
As $x$ is $\left(\frac{1}{4} - \eps \right)n$-balanced, we have $-\left(\frac{1}{2} + 2\eps \right)n \le \deg_R(x) - \deg_B(x) \le \left(\frac{1}{2} + 2\eps \right)n$, and since $c(K_n) = 0$ we get that
\begin{align}
\label{eq:Knminusx}
    \left|c(K_n\setminus\{x\}) \right| = \left|c(K_n) - \deg_R(x) + \deg_B(x) \right| \le \left(\tfrac{1}{2} + 2\eps \right)n.
\end{align}
Then, by the linearity of expectation, we have
\begin{align}
    \left|\mathbb E \left[c(f(F)) \right] \right| &\le \left|\mathbb E\left[c(f(\{v_1\} \cup N(v_1))) \right] + \mathbb E \left[c(f(F\setminus \{v_1\})) \right] \right| \nonumber \\
    &\le \frac{\Delta}{n-1}\cdot\left|\deg_R(x) - \deg_B(x) \right| + \frac{n-\Delta}{\binom{n-1}{2}}\cdot\left|c(K_n\setminus \{x\}) \right| \nonumber \\
    &\le \frac{\Delta}{n-1}\cdot\left|\deg_R(x) - \deg_B(x) \right| + \frac{n(n-\Delta)}{\binom{n-1}{2}}\left(\tfrac{1}{2} + 2\eps \right) \nonumber \\
    &\le \left(\tfrac{1}{2} + 2\eps \right)\Delta + \tfrac{\Delta}{n-1}\left(\tfrac{1}{2} + 2\eps \right) + 1 + 4\eps \nonumber \\
    &\le \tfrac{1}{2}\Delta + 4, \label{eq:small-delta-expectation-inequality}
\end{align}
where the last inequality holds as $n \ge 17$.
Assume, for a contradiction, that there exists no embedding $f : V(F) \rightarrow V(K_n)$ such that $\left|c(f(F)) \right| \le \frac{1}{2}\Delta + 4$.
Then, by inequality \eqref{eq:small-delta-expectation-inequality}, we get that there exist two embeddings $f_B, f_R : V(F) \rightarrow V(K_n)$, both mapping $v_1$ to $x$, such that
\begin{align*}
    c(f_B(F)) < -\tfrac{1}{2}\Delta - 4 && \text{ and } && c(f_R(F)) > \tfrac{1}{2}\Delta + 4.
\end{align*}
Applying \Cref{lem:trivial-algorithm} with $f_B, f_R$ and $\Delta_I \le \Delta'$, we get that there exists a copy $F'$ of $F$ in $K_n$ such that $\left|c(F') \right| < \frac{1}{2}\Delta + 1$, a contradiction.

\paragraph*{Case 2.2: $\Delta' \ge \frac{1}{2}\Delta$.}
Assume without loss of generality that $x$ has red degree $\deg_R(x) \geq \frac{n-1}{2}$.
As $c$ is a balanced colouring, there exists a vertex $y \in V(K_n)$ which has red degree $\deg_R(y) \leq \frac{n-1}{2}$.
We will define an embedding $f: V(F) \rightarrow K_n$ which sends $v_1$ to $x$ and $v_2$ to $y$.

If $\deg_R(y) \ge \frac{1}{4}n$ then we extend $f$ uniformly at random.
Denote $\Delta'' = \deg_F(v_3)$, so we have
\[\Delta'' \le n - 1 - \Delta - \Delta' + 2 \le n + 1 - \tfrac{1}{2}n - \tfrac{1}{4}n \le \tfrac{1}{2}\Delta + 1. \]
Note that, similarly to (\ref{eq:Knminusx}), we also have $|c(K_n \setminus \{ x,y \})| \le 2n$.
Denote $Z = 1_{\{v_1v_2 \in E(F) \}}$.
Then similarly to the analysis in Case \hyperref[case:Delta-prime-small]{2.1}, we get that
\begin{align*}
    \left|\mathbb E\left[c(f(F)) \right] \right| \le\,& 
    \Bigl|\EE\left[c(f(\set{v_1}\union N(v_1))) \right] + 
        \EE\left[c(f(\set{v_2}\union N(v_2))) \right] \nonumber \\
        &+ \EE\left[c(f(F\setminus\{v_1, v_2\})) \right] - \EE\left[c(f(F[\set{v_1,v_2}])) \right] \Bigr| \nonumber \\
    \le\,& 
        \Biggl| \frac{\Delta - Z}{n-2}\Bigl(\deg_R(x) - \deg_B(x) - c(xy)Z \Bigr) \Biggr| +
        \Biggl| \frac{\Delta' - Z}{n-2}\Bigl(\deg_R(y) - \deg_B(y) - c(xy)Z \Bigr) \Biggr| \nonumber \\ &+ 
        \Biggl| \frac{n - \Delta - \Delta' + Z}{\binom{n-2}{2}} \cdot 2n \Biggr| + \bigl| c(xy) \bigr| \nonumber \\ 
    \le\,&
        \frac{\Delta}{n-2} \cdot \max\left(|\deg_R(x) - \deg_B(x)|, |\deg_R(y) - \deg_B(y)| \right) +\frac{\Delta+\Delta'}{n-2} \nonumber \\ &+
        \frac{n - \Delta - \Delta' + 1}{\binom{n-2}{2}} \cdot 2n + 1  \nonumber \\
    \le\,&
        \left(\tfrac{1}{2} + 2\eps \right)\Delta + \frac{2\Delta}{n-2}\left(\tfrac{1}{2}+2\eps \right) +\frac{\Delta+\Delta'}{n-2} + \frac{\frac{1}{4}n + 1}{\binom{n-2}{2}} \cdot 2n+ 1 \nonumber \\
    \le\,&
        \tfrac{1}{2}\Delta + 9 \nonumber,
\end{align*}
where the last inequality follows since for $n \ge 17$ we have
\[2\eps\Delta + \frac{2\Delta}{n-2}\left(\tfrac{1}{2} + 2\eps \right) + \frac{\Delta + \Delta'}{n-2} + \frac{\frac{1}{4}n+1}{\binom{n-2}{2}}\cdot 2n + 1 \le 9. \]
Continuing with an argument similar to the one in Case 2.1, with $\Delta_I \le \Delta'' \le \frac{1}{2}\Delta$ when applying \Cref{lem:trivial-algorithm}, we get that there exists a copy $F'$ of $F$ in $K_n$ with $|c(F')| \le \frac{1}{2}\Delta + 9$.\\

We are left with the case $\deg_R(y) < \frac{1}{4}n$.
In this case we construct an embedding in the following way.
Let $X_R,X_B \subset N_F(v_1)\setminus \{v_2\}$ be two disjoint subsets with $|X_R| = \left\lfloor\frac{3}{8}n \right\rfloor$ and $|X_B| = \left\lfloor\frac{1}{8}n \right\rfloor-1$.
As $F$ is a forest we have $\left|N_F(v_2) \cap \left(N_F(v_1) \cup \{v_1\} \right)\right| \le 1$.
Moreover, we have $\deg_F(v_2) \ge \frac{1}{4}n$, so we may choose a subset $Y_B \subset N_F(v_2) \setminus \left(N_F(v_1) \cup \{v_1\} \right)$ with $|Y_B| = \left\lfloor\frac{1}{4}n \right\rfloor-1$.
Since $\deg_B(y) > \frac{3}{4}n$, for any $f : V(F) \rightarrow K_n$ we have that $\left|N_B(y) \setminus (f(X_R)\cup f(X_B)) \right| \ge \frac{1}{4}n$.
Recall that we require $f$ to map $v_1$ to $x$ and $v_2$ to $y$.
Hence, since $x$ is $\left(\frac{1}{4}-\varepsilon \right)n$-balanced we may let $f$ be an embedding of $F$ such that
\begin{itemize}
    \item $f(X_R) \subset N_R(x)$,
    \item $f(X_B) \subset N_B(x)$, and
    \item $f(Y_B) \subset N_B(y) \setminus \left(f(X_R) \cup f(X_B) \right)$.
\end{itemize}
We then get an embedding $f$ of $F$ with at least $\left\lfloor\frac{3}{8}n \right\rfloor$ red edges and at least $\left\lfloor\frac{3}{8}n \right\rfloor-2$ blue edges, so that $|c(f(F))| \le \frac{1}{4}n + 4 \le \frac{1}{2}\Delta + 4$.
\end{proof}

\Cref{thm:spanning-forest-general} follows from \Cref{thm:spanning-forest} as a simple corollary.
\begin{proof}[Proof of \Cref{thm:spanning-forest-general}]
    It suffices to consider $15 < \Delta < \frac{1}{2}n$.
    It is easy to show by simple algebraic manipulations that for any $|x| \le \frac{1}{4}n$, we have
    \[\sqrt{\left(\tfrac{x}{2} \right)^2+4n} \le \tfrac{1}{8}n + 16. \]
    Applying this inequality to $x=\Delta - \tfrac{1}{4}n + 1$, we obtain
    \begin{align*}
        \tfrac{1}{2}\left(\Delta - \tfrac{1}{4}n + 3 \right) + \sqrt{\tfrac{1}{4} \left(\Delta - \tfrac{1}{4}n + 1 \right)^2 + 4n} &\leq \tfrac{1}{2}\left(\Delta - \tfrac{1}{4}n + 3 \right) + \tfrac{1}{8}n+16 \\
        &\leq \tfrac{1}{2}\Delta + 18,
    \end{align*}
    which, together with \Cref{thm:spanning-forest}, implies (\ref{eq:conjUB}).
\end{proof}

We end this section with a proof of \Cref{cor:main} as a direct consequence of \Cref{thm:spanning-forest}.
\begin{proof}[Proof of \Cref{cor:main}]
    Let $\xi \in \left[-\frac{1}{4} + \frac{16}{n}, -\eta \right]$ be such that $\xi n = \Delta - \frac{1}{4}n$.
    By \Cref{thm:spanning-forest} and Bernoulli's inequality, we get that there exists a copy $F'$ of $F$ in $K_n$, for $n$ sufficiently large, such that
    \begin{align*}
        \left|c(F') \right| &\le \frac{2+\xi n}{2} + \frac{\sqrt{(2+\xi n)^2 + 16n}}{2} \\
        &= \frac{2+\xi n}{2} - \frac{2+\xi n}{2} \sqrt{1 + \frac{16n}{(2+\xi n)^2}} \\
        &\le \frac{2+\xi n}{2} - \frac{2+\xi n}{2} \left(1 + \frac{8n}{(2+\xi n)^2} \right) \\
        &= \frac{4n}{|2+\xi n|} \\
        &\le \frac{4n}{\eta n-2} \\
        &\le (1+o(1))\frac{4}{\eta},
    \end{align*}
    as wanted.
\end{proof}


\section{Concluding remarks and further work}
\label{sec:conclusion}

The first natural further line of research would be to obtain tight bounds concerning \Cref{thm:spanning-forest} for $\Delta \ge \frac{1}{4}n$.
If $n$ is a multiple of $4$, the following construction was noted in \cite{mohr2022zero}.
We let $V(K_n)= \{u_i: i \in I \} \cup \{v_i: i \in I \}$ for $I=[n/2]$ and we let $c_0\from E(K_n) \to \set{-1, +1}$ be defined as
\begin{align}
\label{eq:ConstructionMPR}
    c_0(xy)=\begin{cases}
        -1 &\text{if } xy=u_iu_j \text{ for some } i,j \in I, \text{ or } xy=u_iv_j \text{ for some } i,j \in I \text{ s.t.\ } i+j \text{ is odd, } \\
        +1 &\text{otherwise.}
    \end{cases}
\end{align}
If $F$ is an $n$-vertex star, then every embedding has sum precisely equal to $\frac{\Delta-1}{2}$ in absolute value.
Furthermore, if $F$ is close to a star -- that is $F$ has one vertex of degree $\Delta\ge \frac{3}{4}n$ and some other arbitrary vertices and edges -- then the absolute value of the sum of any embedding in colouring $c_0$ is at least $2\left(\Delta - \frac{3}{4}n \right) + 1$.
It is worth noting that these lower bounds are far from our upper bounds. 
In particular, we have no nontrivial lower bound when $\frac{1}{4}n \le \Delta < \frac{3}{4}n$. \\

One different line of further research would be to determine whether a form of \Cref{thm:spanning-forest} still holds when the complete graph is replaced by a dense enough host graph.
Indeed, Koml{\'o}s, S{\'a}rk{\"o}zy, and Szemer{\'e}di \cite{KSS01} have proved that if $\delta>0$, and $G$ is a graph on $n$ vertices and minimum degree $\delta(G) \ge \left(\frac{1}{2}+\delta \right)n$, then $G$ contains every spanning tree $T$ with maximum degree $\Delta(T) = O\left(\frac{n}{\log n} \right)$, provided that $n$ is sufficiently large. 
Whether a low-sum copy of $T$ can be found in a graph $G$ with high minimum degree equipped with a balanced $2$-colouring is an open question. \\

It would also be interesting to see whether an analogue of the main result can still be obtained when the complete graph is coloured with a small imbalance.
We remark that in this case the bound would depend on the perturbation of the $2$-colouring, as our final result shows.

\begin{theorem}
\label{thm:UBZeroSumPerturbation}
    For every sufficiently small $\varepsilon >0$ and sufficiently large $n$ there exists a red-blue colouring of the edges of $K_n$ such that the density of the red edges is $\frac{1}{2}-\varepsilon \le d_R \le \frac{1}{2}+\varepsilon$ and every copy $T$ of a star on $n$ vertices in $K_n$ satisfies $|c(T)| \geq (\frac{1}{2}+\varepsilon^2)n \geq (\frac{1}{2}+\varepsilon^2) \Delta$.
\end{theorem}

We prove \Cref{thm:UBZeroSumPerturbation} by using a modification of the colouring \eqref{eq:ConstructionMPR} in \Cref{appendixC}. \\

One final line of research would be to investigate whether a similar bound as \Cref{thm:spanning-forest-general} could be obtained if we do not restrict the graph $F$ to be a forest.
Indeed, the argument by Mohr, Pardey, and Rautenbach \cite{mohr2022zero} mentioned in the introduction and slightly generalised in \Cref{lem:trivial-algorithm} show that, for every graph $H$ on $n$ vertices, maximum degree $\Delta_H$ and minimum degree $\delta_H$, and every balanced colouring of $K_n$, then there exists a copy $H'$ of $H$ in $K_n$ such that $|c(H')| \leq \Delta_H+\delta_H \leq 2\Delta_H$.
It would be interesting to see whether this simple argument is in some sense optimal, or whether tighter bounds could be obtained.
More precisely, we ask the following.
\begin{question}
\label{qu:Non-forest}
    Let $H$ be a graph on $n$ vertices, with maximum degree $\Delta=\Delta(H)$, and let $c : E(K_n) \rightarrow \{-1, +1\}$ be a balanced colouring of the edges of $K_n$. 
    Does there exists a copy $H'$ of $H$ in $K_n$ such that
    \begin{align}
    \label{eq:quest1}
        \left|c(F') \right| \le \alpha \Delta,
    \end{align}
    for some universal constant $\alpha < 2$?
\end{question}

\paragraph*{Acknowledgements.}
The authors would like to thank their PhD supervisor Professor B\'{e}la Bollob\'{a}s for his support and valuable comments, and for a stimulating discussion which led to \Cref{qu:Non-forest}.
The first and second authors are funded by Trinity College, Cambridge.
The third author is funded by EPSRC (Engineering and Physical Sciences Research Council) and by the Cambridge Commonwealth, European and International Trust.

\bibliographystyle{abbrvnat}  
\renewcommand{\bibname}{Bibliography}
\bibliography{main}


\appendix

\section{Appendix: Proof of Lemma \ref{lem:ExistenceCbalancedVertices}}
\label{appendixA}

\begin{proof}[Proof of \Cref{lem:ExistenceCbalancedVertices}]
Set $V \coloneqq V(K_n)$, and define the following three subsets of $V$:
\begin{align*}
    D^{-} &\coloneqq \left\{ v \in V ~:~ \deg_R(v) < \left(\tfrac{1}{4}-\eps\right)n \right\}, \\
    D^{+} &\coloneqq \left\{v \in V ~:~ \deg_R(v) > \left(\tfrac{3}{4}+\eps \right)n \right\}, \\
    D &\coloneqq \left\{v \in V ~:~ \left(\tfrac{1}{4}-\eps \right)n \le \deg_R(v) \le \left(\tfrac{3}{4}+\eps \right)n \right\}.
\end{align*}
Note that $D^- \cup D^+ \cup D$ is a partition of the vertices of $K_n$.
By symmetry, we may assume that $|D^{+}| \le |D^{-}|$, and consequently that $|D^{-}| \ge \frac{1}{2}(n-|D|)$.
Then the number of red edges in the graph satisfies
\begin{align*}
   \tfrac{1}{2}\tbinom{n}{2}=e_R &= e_R(D^-,V)+e_R(V\setminus D^-, V\setminus D^-)\\
   &\le \left(\tfrac{1}{4} - \eps \right)n|D^{-}|+\tbinom{n-|D^-|}{2}.
\end{align*}
Note that the function $g(u) = \left(\frac{1}{4}-\eps \right)nu + \binom{n-u}{2}$ is a quadratic polynomial in $u$ with a positive leading coefficient, and therefore attains its maximum on an interval $[a,b]$ at either $a$ or $b$.
As $\frac{1}{2}(n-|D|) \le |D^-| \le n$, it follows that $g(|D^-|) \le \max\left(g\left(\frac{1}{2}(n-|D|) \right), g(n) \right)$.
It can be easily verified that for any $n \ge 5$, since $|D|\ge 0$, we have
\begin{align*}
    g\left(\tfrac{1}{2}(n-|D|) \right) = \left(\frac{1}{4} - \eps \right)\frac{n(n-|D|)}{2} + \binom{\tfrac{1}{2}(n+|D|)}{2} \ge \left(\frac{1}{4} - \eps \right)n^2 = g(n),
\end{align*}
and hence
\begin{align*}
    \frac{1}{2} \binom{n}{2} = e_R \le g(|D^-|) \le g\left(\tfrac{1}{2}(n-|D|) \right) = \left(\frac{1}{4} - \eps \right) \frac{n(n-|D|)}{2} + \binom{\tfrac{1}{2}(n+|D|)}{2}.
\end{align*}
Rearranging, we get 
\begin{align*}
    \frac{1}{4}|D|^2 + \left(\left(\frac{1}{4} + \eps \right)n - \frac{1}{2} \right)|D| - \eps n^2 \ge 0,
\end{align*}
which implies
\begin{align*}
    2\abs{D}\geq \sqrt{(1 + 24\eps + 16\eps^2)n^2 - 4(1 + 4\eps)n + 4} - (1 + 4\eps)n + 2.
\end{align*}
Thus to show that $\abs{D}\geq \eps n + 1$, it suffices to show that 
$$(1 + 24\eps + 16\eps^2)n^2 - 4(1 + 4\eps)n + 4 \geq (1 + 6\eps)^2 n^2.$$
Rearranging, it suffices to show that 
\begin{align*}
    (3\eps - 5\eps^2)n^2 - (1 + 4\eps)n + 1 \geq 0.
\end{align*}
As $0 < \eps < \frac{1}{4}$, we have $(3\eps - 5\eps^2) > 0$, and therefore it suffices to show that
\begin{align}
\label{eq:BoundSufficientn}
    n \geq \frac{1+4\eps+\sqrt{1-4\eps+36\eps^2}}{2\eps(3-5\eps)}.
\end{align}
It is easy to show that $\frac{1+4\eps+\sqrt{1-4\eps+36\eps^2}}{2(3-5\eps)} \le 1$ for every $0< \eps < \frac{1}{4}$.
Therefore our assumption $n \ge \frac{1}{\eps}$ implies that \eqref{eq:BoundSufficientn} is indeed satisfied, finishing the proof.
\end{proof}

\section{Appendix: Proof of Lemma \ref{lem:optim}}
\label{appendixB}

\begin{proof}[Proof of \Cref{lem:optim}]
    Let $f(\eps)=\frac{2}{\eps}+1$ and $g_{\xi,n}(\eps)=2+(\xi n + 2\eps n)_+$ so that 
    $$\phi_n(\xi)= \inf_{\frac{1}{n} \le \eps \le \frac{1}{8}}\max\left(f(\eps), g_{\xi,n}(\eps) \right).$$
    It is clear that $f$ is decreasing in $\varepsilon$, and $g_{\xi,n}$ is non-decreasing.
    Moreover, since $\xi \ge -\frac{1}{4} + \frac{16}{n}$, we have
    \begin{align*}
        (f - g_{\xi,n})\left(\tfrac{1}{n} \right) > 0 && \text{and} && (f - g_{\xi,n})\left(\tfrac{1}{8} \right) < 0.
    \end{align*}
    Hence, as both $f$ and $g_{\xi,n}$ are continuous for $\eps>0$, there exists $ \eps_*$ such that $\frac{1}{n} \le \eps_* \le \frac{1}{8}$ and $f(\eps_*) = g_{\xi,n}(\eps_*) = \phi_n(\xi)$.
    In fact, since $\frac{2}{\eps_*}+1 > 2$, we get that $\eps_*$ satisfies
    \[2 + \xi n + 2\eps_* n = \frac{2}{\eps_*} + 1. \]
    Rearranging this and solving the second degree equation in $\eps_{*}$, we obtain
    \begin{align*}
        \eps_{*}=\frac{-1-\xi n+\sqrt{(1+\xi n)^2+16n}}{4n},
    \end{align*}
    implying
    \begin{align}
        \phi_n(\xi) = 2 + \xi n + 2\varepsilon_{*}n = \frac{3+\xi n}{2} + \frac{\sqrt{(1+\xi n)^2+16n}}{2},
    \end{align}
    as required.
\end{proof}

\section{Appendix: Proof of Theorem \ref{thm:UBZeroSumPerturbation}}
\label{appendixC}

\begin{proof}[Proof of \Cref{thm:UBZeroSumPerturbation}]
    It suffices to prove that for $n$ sufficiently large there exists a red-blue
    colouring of the edges of $K_n$ such that the density of the red edges is $\frac{1}{2}-\varepsilon \le d_R \le \frac{1}{2}+\varepsilon$ and every vertex $v$ satisfies $\deg_R(v) \ge \left(\frac{3}{4}+\frac{1}{2}\varepsilon^2 \right)n$ or $\deg_R(v) \le \left(\frac{1}{4}-\frac{1}{2}\varepsilon^2 \right)n$.
    For this purpose, we generalise the construction \eqref{eq:ConstructionMPR} given in \cite{mohr2022zero}.
    We let $A$ and $B$ be a partition of $[n]$ such that $|A| = \left(\frac{1}{2}-\varepsilon \right)n$ and $|B| = \left(\frac{1}{2}+\varepsilon \right)n$.
    We colour blue those edges with both endpoints in $A$ and red those edges with both endpoints in $B$.
    Let $d = \frac{x}{y}$ be a rational number such that $0 < d < 1$ to be chosen later.
    Writing $A = \{ a_1, \dots, a_{|A|} \}$ and $B= \{ b_1, \dots, b_{|B|} \}$ we define
    \begin{align*}
    c(a_ib_j)=\begin{cases}
        \text{red} &\text{ if } \text{mod}_y(i+j) \leq x, \\
        \text{blue} &\text{ if } \text{mod}_y(i+j) > x,
    \end{cases}
\end{align*}

    where by $\text{mod}_y(i+j)$, we mean the number in the set $[y]$ which is congruent to $i+j$ modulo $y$.
    Note that the desired conditions on red/blue degrees can be rewritten as
    \begin{align}
    \label{eq:RedDegree}
        \left(\tfrac{1}{2} + \varepsilon \right) + d\left(\tfrac{1}{2} - \varepsilon \right) \ge \tfrac{3}{4} + \tfrac{1}{2}\varepsilon^2 + O_d \left(\tfrac{1}{n}\right)
    \end{align}
    for the red degree, and
    \begin{align}
    \label{eq:BlueDegree}
        \left(\tfrac{1}{2} - \varepsilon \right)+(1-d)\left(\tfrac{1}{2} + \varepsilon \right) \ge \tfrac{3}{4} +\tfrac{1}{2}\varepsilon^2 + O_d \left(\tfrac{1}{n}\right)
    \end{align}
    for the blue degree.
    Moreover, note that the density of red edges in this graph is
    $$d_R = 2d\left(\tfrac{1}{2} - \varepsilon \right) \left(\tfrac{1}{2} + \varepsilon \right) + \left(\tfrac{1}{2} +\varepsilon \right)^2 + O_d \left(\tfrac{1}{n}\right),$$
    so the conditions on the density of red edges can be rewritten as
    \begin{align}
    \label{eq:EdgesDensity}
        \tfrac{1}{4} - \tfrac{1}{2}\varepsilon \le d\left(\tfrac{1}{2} - \varepsilon \right)\left(\tfrac{1}{2} + \varepsilon \right) + \tfrac{1}{2}\left(\tfrac{1}{2} + \eps \right)^2 + O_d \left(\tfrac{1}{n}\right) \le \tfrac{1}{4} + \tfrac{1}{2}\varepsilon.
    \end{align}
    Thus, to prove our theorem, it suffices to show that there is a rational number $0 < d < 1$ which simultaneously satisfies \eqref{eq:RedDegree}, \eqref{eq:BlueDegree} and \eqref{eq:EdgesDensity}.
    By \eqref{eq:RedDegree}, \eqref{eq:BlueDegree} we get that such $d$ should satisfy
    \begin{align}
    \label{eq:d-interval-1}
        \frac{1/4+\varepsilon^2/2-\varepsilon}{1/2-\varepsilon} < d < \frac{1/4-\varepsilon^2/2}{1/2+\varepsilon},
    \end{align}
    and by \eqref{eq:EdgesDensity} such $d$ should satisfy
    \begin{align}
    \label{eq:d-interval-2}
        \frac{1/8-\varepsilon-\varepsilon^2/2}{1/4-\varepsilon^2} < d < \frac{1}{2}.
    \end{align}
    To show the existence of such $d$, it suffices to check that the two open intervals of values of $d$ given by \eqref{eq:d-interval-1} and \eqref{eq:d-interval-2} are non empty and moreover, have a non empty intersection.
    In fact, it is easy to check that for any $0 < \eps < \frac{1}{2}$ we have
    \begin{align*}
        \frac{1/8-\varepsilon-\varepsilon^2/2}{1/4-\varepsilon^2} < \frac{1/4+\varepsilon^2/2-\varepsilon}{1/2-\varepsilon} < \frac{1/4-\varepsilon^2/2}{1/2+\varepsilon} < \frac{1}{2},
    \end{align*}
    meaning that the open interval given by \eqref{eq:d-interval-1} is non empty and is completely contained in the interval given by \eqref{eq:d-interval-2}.
    This finishes the proof.
\end{proof}

\end{document}